\documentstyle[12pt]{article}

\advance\textheight by -\headheight
\advance\textheight by -\headsep

\def\_#1{\mathop{\hspace{-2pt}^{}_{#1}}}
\def\suml  {\mathop{\sum}  \limits}

\def\X{\mathop{\cal X}\nolimits}

\def\R{\mathop{\cal R}\nolimits}
\def\P{\mathop{\cal P}\nolimits}
\def\W{\mathop{\cal W}\nolimits}
\def\L{\mathop{\cal L}\nolimits}
\def\sign{\mathop{\rm sign}\nolimits}
\def\smp{\mathop{\scriptscriptstyle (+)}\nolimits}
\def\smm{\mathop{\scriptscriptstyle (-)}\nolimits}
\def\ve{\mathop{\varepsilon}}

\def\QA{\mathop{\rm QA}\nolimits}
\def\WQA{\mathop{\rm WQA}\nolimits}
\def\cdc{,\ldots,}
\def\sameauthor{\rule[.27em]{6.2ex}{.4pt}\ }
\def\r{\mathop{\rho}\nolimits}
\def\gro  {\mathop{\succ_\rho}\nolimits}
\def\gero {\mathop{\succeq_\rho}\nolimits}
\def\grop {\mathop{\succ_{\rho'}}\nolimits}

\def\sro{\mathop{\sim_\rho}\nolimits}
\def\srop{\mathop{\sim_{\rho'}}\nolimits}

\def\flushpar{\par\noindent}
\def\lfi{\flushpar\leftskip=1.25em\hspace{-1.25em}}

\begin{document}

\centerline{\huge From Incomplete Preferences to Ranking}
\centerline{\huge via Optimization}

\bigskip\bigskip\medskip

\centerline{\Large Pavel~Chebotarev\footnote{Corresponding author. E-mail: \tt chv@lpi.ru} and Elena Shamis}
\begin{center}
Institute of Control Sciences of the Russian Academy of Sciences\\
65 Profsoyuznaya Str., Moscow 117997, Russia
\end{center}

\bigskip\bigskip\medskip

\vspace{5pt}
{\small {\bf Abstract:} We consider methods for aggregating preferences
that are based on the resolution of discrete optimization problems. The
preferences are represented by arbitrary binary relations (possibly
weighted) or incomplete paired comparison matrices. This incomplete
case remains practically unexplored so far. We examine the properties
of several known methods and propose one new method. In particular, we
test whether these methods obey a new axiom referred to here as {\it
Self-Consistent Monotonicity.} Some results are established that
characterize solutions of the related optimization problems.
}
\bigskip

\flushpar{\bf Keywords:} Aggregation of preferences, ranking, paired
comparisons, quadratic assignment problem, Kemeny median
\bigskip

\section{Introduction}

We consider methods for aggregating preferences that are based on the
resolution of discrete optimization problems. For a review and
references see Cook and Kress (1992), and Belkin and Levin (1990), and
also David (1988) and Van Blokland-Vogelesang (1991). Some algorithmic
aspects can be found in Barth\'elemy (1989) and Litvak (1982). The
preferences are represented by arbitrary binary relations (possibly
weighted) or incomplete paired comparison matrices.  The outcome of an
aggregation method is a set of ``optimal'' rankings (linear or weak
orders) of the alternatives.  Namely, a ranking is said to be optimal if it
provides an extremum of some chosen objective function that expresses
the connection (or proximity) between an arbitrary ranking and the
original preferences. One special feature of the aggregation problem
with incomplete preferences is that the Borda-like score is not
any more a good rating index, since it does not take into account the
number of comparisons and the strength of ``opponents" of each alternative.
The incomplete case remains practically unexplored so far. In this
paper, we present some initial results concerning this problem.

We examine the properties of several known methods formulated for the
incomplete case and propose one new method. In particular, we test
whether these methods obey a new axiom, referred to as {\it
Self-Consistent Monotonicity.} Our results suggest that the methods
under consideration hardly satisfy this condition because of their
discreteness. This paper provides only one example of methods
satisfying Self-Consistent Monotonicity (Section~5). Actually, that
method is an {\it indirect scoring procedure} rather than an
aggregating operator based on discrete optimization. Indirect scoring
procedures are considered in Chebotarev and Shamis (1996), where a
sufficient condition of Self-Consistent Monotonicity and some more
positive examples are given. A discussion of some other axioms for
aggregating incomplete preferences can be found in Chebotarev and
Shamis (1994).

\section{An illustrative example and the aims of the paper}

We start with a simple example. Let there be four candidates and
two voters. The preferences of these voters are incomplete.
Namely, the first voter says: ``$X_2$ and $X_3$ are better than
$X_4$, and I know nothing about $X_1$." (Fig.\ 1 a.) The second
voter says: ``$X_1$ is better than $X_3$, and I know nothing about
$X_2$ and $X_4$." (Fig.\ 1 b.)

\begin{figure}[htb]
\setlength{\unitlength}{5mm}
\begin{picture}(26,15.6)
\put(2,4){\line(1,0){7}}
\put(2,14){\line(1,0){7}}
\put(2,4){\line(0,1){10}}
\put(9,4){\line(0,1){10}}
\put(2,4){\begin{picture}(7,10)
             \put(1.5,8.0){\circle{0.2}} 
             \put(1.5,4.5){\circle{0.2}} 
             \put(5.5,7.0){\circle{0.2}} 
             \put(5.5,2.5){\circle{0.2}} 
             \put(1.2,8.3){$X_1$}
             \put(1.1,3.4){$X_3$}
             \put(5.2,7.3){$X_2$}
             \put(5.2,1.6){$X_4$}
             \put(1.8,4.37){\vector(2,-1){3.3}} 
             \put(5.5,6.7){\vector(0,-1){3.7}}  
             \put(0.62,-1){\small a. The first voter}
          \end{picture}}
\put(11,4){\line(1,0){7}}
\put(11,14){\line(1,0){7}}
\put(11,4){\line(0,1){10}}
\put(18,4){\line(0,1){10}}
\put(11,4){\begin{picture}(7,10)
             \put(1.5,8.0){\circle{0.2}}
             \put(1.5,4.5){\circle{0.2}}
             \put(5.5,7.0){\circle{0.2}}
             \put(5.5,2.5){\circle{0.2}}
             \put(1.2,8.3){$X_1$}
             \put(1.1,3.4){$X_3$}
             \put(5.2,7.3){$X_2$}
             \put(5.2,1.6){$X_4$}
             \put(1.5,7.7){\vector(0,-1){2.8}} 
             \put(0.1,-1){\small b. The second voter}
          \end{picture}}
\put(20,4){\line(1,0){7}}
\put(20,14){\line(1,0){7}}
\put(20,4){\line(0,1){10}}
\put(27,4){\line(0,1){10}}
\put(20,4){\begin{picture}(7,10)
             \put(1.5,8.0){\circle{0.2}}
             \put(1.5,4.5){\circle{0.2}}
             \put(5.5,7.0){\circle{0.2}}
             \put(5.5,2.5){\circle{0.2}}
             \put(1.2,8.3){$X_1$}
             \put(1.1,3.4){$X_3$}
             \put(5.2,7.3){$X_2$}
             \put(5.2,1.6){$X_4$}
             \put(1.8,4.37){\vector(2,-1){3.3}} 
             \put(1.5,7.7){\vector(0,-1){2.8}}  
             \put(5.5,6.7){\vector(0,-1){3.7}}  
             \put(1.2,  -1){\small c. combined}
             \put(1.2,-1.7){\small~~(two voters)}
          \end{picture}}
\put(8.2,1.0){\small{\sc Figure 1.}  Preferences of two voters}
\end{picture}
\end{figure}

Certainly, these preferences are extremely poor, and it is
difficult to make a decision based on them. Nevertheless, having
in mind more lifelike situations, we may pose a question: Which
principles should be followed when aggregating incomplete
preferences. Every answer will imply some consequences applicable
to complete preferences as well.

Obviously, in this example $X_1$ or $X_2$ should be the first. We
believe that $X_1$ has a small advantage over $X_2$, since $X_1$
defeats a ``stronger" opponent. The rank order of the other
candidates is more clear: $X_2\succ X_3\succ X_4.$

Our main requirement to the aggregating operators comes down to
the following. Suppose we consider two alternatives, and the
first alternative as compared to the second one
\newline
-- achieves the better scores against the ``stronger" opponents
or
\newline
-- achieves the better scores against the opponents of the same
``strengths" or
\newline
-- achieves the same scores against the ``stronger" opponents.

Then the first alternative should be placed {\it higher} than the
second one in the social ranking.

Now we have to explain what is meant by ``stronger". In the
above requirement, ``stronger" signifies ``{\it is placed higher
in the social ranking that is mentioned in the requirement}."
This requirement is formalized in Section~7 and is referred to as
{\it Self-Consistency}. To obtain {\it Self-Consistent
Monotonicity}, the main axiom in this paper, we require that if the
first alternative additionally achieves some extra ``wins" and/or
the second alternative has extra ``losses", then the first
alternative should remain {\it higher} than the second one in the
social ranking.

It turns out that very few methods satisfy these natural axioms. In
this paper, we adjust a number of discrete optimization procedures to
the case of incomplete preferences. The most familiar are the close
procedures by Kemeny (1959) and Slater (1961) (their idea has been
initially suggested by Kendall (1955)).  For instance, the Slater
method minimizes the number of arcs (which designate individual
binary preferences) directed upwards (from a ``worse" alternative to a
``better" one) in the social ranking.  When applied to the above
example, this method produces three optimal social rankings (with no
upward arcs), namely, $X_1\succ X_2\succ X_3\succ X_4$, $X_2\succ
X_1\succ X_3\succ X_4$, and $X_1\succ X_3\succ X_2\succ X_4$.  Only
the first of them preserves Self-Consistent Monotonicity.

The aims of this paper are:
\newline
-- to collect some discrete optimization methods, to represent
them in the assign\-ment-like form (Section~4), and to define
their modifications that generate weak orders (Section~9);
\newline
-- to give some results characterizing solutions of these
optimization problems (Section~6);
\newline
-- to introduce Self-Consistent Monotonicity (Section~7);
\newline
-- to present some necessary conditions of Self-Consistent
Monotonicity (Sections 8 and 10);
\newline
-- to prove Self-Consistent Monotonicity for the generalized row
sum method and to test discrete optimization methods in this
respect (Sections 8, 10, and 11);
\newline
-- to outline the difficulties in obeying Self-Consistent
Monotonicity by discrete optimization methods (Section~11).

\section{Notation}

The case of incomplete and possibly weighted preferences requires
some more complex notation. Let $\X=\{X_1\cdc X_n\}$ be a set of
alternatives to be compared. To represent arbitrary preference
relations of individuals, both ordinary binary relations and weighted
ones, we use {\it incomplete paired comparison matrices.} Such a
matrix of the $p$th individual $(p=1\cdc m)$ is an $n\times n$ table
$R^{(p)}=(r_{ij}^p)_{i,j=1}^n$ whose entries $r_{ij}^p$ and
$r_{ji}^p$ represent the result of comparing $X_i$ to $X_j$ by this
individual. The ordered pair $(r_{ij}^p, r_{ji}^p)$ will be called
the {\it outcome\/} of that comparison. Here the value of $r_{ij}^p$
can be interpreted as a measure of advantage of $X_i$ over $X_j$
(something like the number of scored goals in sport). If $X_i$ and
$X_j$ have not been compared by that individual, the two
corresponding cells of the table remain empty (i.e., $r_{ij}^p$ and
$r_{ji}^p$ are undefined). The diagonal elements $r_{ii}^p\; (i=1\cdc
n)$ do not correspond to any comparisons and are defined according to
some convention. The collection of matrices $\R=(R^{(1)}\cdc
R^{(m)})$ is called {\it an array of paired comparisons} of the
alternatives $X_1\cdc X_n$.  If all $R^{(p)}\;(p=1\cdc m)$ are
completely defined, then $\R$ is said to be {\it complete}.

Let $R=(r_{ij})$ be a complete $n\times n$-matrix, where
$$
r_{ij}=\suml_{p|(i,j)}r_{ij}^p,\quad i,j=1\cdc n
$$
with summation over $p$ for which $r_{ij}^p$ is defined in $\R$.
If for every $p$, $r_{ij}^p$ is undefined in $R^{(p)}$, then, by
definition, $r_{ij}=0$. The total number of comparisons between
$X_i$ and $X_j\ne X_i$ will be denoted by $m_{ij}$:
$$
m_{ij}=\bigl\vert\{p: r_{ij}^p\;\; {\em is\;\; defined\;\; in\;}
\R\}\bigr\vert.
$$
If $i=j$ then $m_{ij}=0.$

An {\it aggregating operator\/} is a mapping assigning to every $\R$
with fixed $n\ge2$ and $m\ge1$ a nonempty set of weak orders on $\X$.
These weak orders are called {\it optimal}. If an aggregating
operator always generates only {\it linear\/} optimal orders, we say
that it is {\it strict}. Recall that a {\it weak order\/} is a
complete and transitive binary relation; a {\it linear order\/} is an
antisymmetric weak order.  $\W$ and $\L$ will denote the set of all
weak orders on $\X$ and the set of all linear orders on $\X$,
respectively. These binary relations are considered here as those of
preference, i.e., for such a relation $\r$, $X_i\r X_j$ means ``$X_i$
is not worse than $X_j$ according to $\r$."

For any binary relation $\r$ on $\X$ and for every $X_i\in\X$, the
{\it Copeland index of $X_i$ in\/} $\r$ can be defined as follows:
$$
\r(X_i)=\r(i)=\r i=\bigl\vert\{X_j\in\X:\: X_i\r X_j\}\bigr\vert-
\bigl\vert\{X_j\in\X:\: X_j\r X_i\}\bigr\vert.
\eqno(1)
$$

If $\r$ is a weak order, then $\r i>\r j$ can be interpreted as ``$X_i$
is better than $X_j$ in $\r$". In this case we write $X_i\gro X_j$. If
$\r$ is a weak order and $\r i=\r j$, we say that $X_i$ and $X_j$ are
{\it tied\/} in $\r$ and write $X_i\sro X_j$. The expression $X_i\gero
X_j$ denotes the disjunction of $X_i\gro X_j$ and $X_i\sro X_j$. Thus
the Copeland index enables to extend the notion of rank to weak
orders (in the manner like that used in statistics).

Paired comparisons can be dichotomous ($r_{ij}^p\in\{-1,1\}$ or
$r_{ij}^p\in\{0,1\}$), with draws ($r_{ij}^p\in\{-1,0,1\}$ or
$r_{ij}^p\in\{0,{1\over2},1\}$), numerical, and so on; different
connections between $r_{ij}^p$ and $r_{ji}^p$ can be imposed. In this
paper, we consider incomplete paired comparisons with the only
connection between $r_{ij}^p$ and $r_{ji}^p$ that if $r_{ij}^p$ is
defined in $\R$, then $r_{ji}^p$ is defined too. Let us suppose that
there exist $r_{\min}$ and $r_{\max}>r_{\min}$ such that all entries
$r_{ij}^p$ must belong to the closed interval $[r_{\min},r_{\max}]$.
An outcome $(r_{ij}^p, r_{ji}^p)$ of comparing $X_i$ to $X_j$ will be
called a {\it maximal win} if $r_{ij}^p=r_{\max}$ and
$r_{ji}^p=r_{\min};$ it will be called a {\it maximal loss\/} if
$r_{ij}^p=r_{\min}$ and $r_{ji}^p=r_{\max}$. As it has been mentioned
above, we do not require that $\R$ consists of only maximal wins and
maximal losses, however our results are applicable to the paired
comparisons of that type as well. We only suppose that maximal wins
(maximal losses) {\it are admissible.} By definition, put
$r_{ii}^p=0$ for all $i\in\{1\cdc n\}$ and $p\in\{1\cdc m\}$.

Now we introduce the {\it Copeland index $t(X_i)$ of $X_i$ in the array
of paired comparisons\/} $\R$:
$$
t(X_i)=t_i=\suml_{(j,p)|i}(r_{ij}^p - r_{ji}^p)=\suml_{j=1}^n
(r_{ij}-r_{ji}),
\eqno(2)
$$
where $(j,p)|i$ denotes summation over $j$ and $p$, for which
$r_{ij}^p$ is defined in $\R$.

\medskip
{\it Remark.} This framework tolerates many diverse ways of extracting
the numbers $r_{ij}^p$ from the individual perceptions. In this paper,
we confine ourselves to the data for which sums and differences (such
as in the above formula) make sense. This means that they are
compatible with the scale type or are meaningful in some other exact
model. Below we examine the properties of different objective
functions based on these operations.

\section{Objective functions for aggregating\newline preferences}

It can be easily shown that many known optimization methods for
aggregating preferences can be reduced to {\it quadratic assignment
problems\/} of the form
$$
\QA(R,C):\;\; {\em maximize} \;\;\suml_{i=1}^n\suml_{j=1}^n r_{ij}C(\r
i,\r j)\;\;{\em for}\; \r\in\P,
\eqno(3)
$$
where $C(\cdot,\cdot)$ is a fixed {\it structure function}, and $\r
i$ is the Copeland index of $X_i$ in $\r$. Now $\P$ is the set of
all {\it linear orders\/} on $\X$ (i.e., $\P=\L$), but below
a more general case ($\P=\W$) is considered too. The quadratic
assignment objective function measures some multiplicative
consistency (depending on $C(\cdot,\cdot)$) between the original
preferences and a tentative resulting order $\r$.

The formulation (3) of the quadratic assignment problem is not
conventional. We use the Copeland index of $X_i$ in $\r$ instead of a
simple rank (see, e.g., Hubert (1976) and Arditti (1983)) since this
straightforward generalization provides an easy way to introduce {\it
weak quadratic assignment problems\/} ($\WQA(R,C)$) involving
arbitrary {\it weak orders\/} ($\P=\W$) instead of {\it linear
orders\/} on $\X$.

As long as quadratic assignment problems with $\P=\L$, are considered,
the following structure functions are relevant:
$$
C_1(x,y)=\sign(x-y),\quad
C_2(x,y)=(\sign(x-y))^{\smp},\quad
C_3(x,y)=(\sign(x-y))^{\smm},\quad
$$
$$
C_4(x,y)=x-y,\quad
C_5(x,y)=(x-y)^{\smp},\quad
C_6(x,y)=(x-y)^{\smm},\quad
$$
where $z^{\smp}=\max(z,0)\:$, $z^{\smm}=\min(z,0)$, and
$$
\sign z=\cases{         -1, &if $z<0$,\cr
              \phantom{-}0, &if $z=0$,\cr
              \phantom{-}1, &if $z>0$.\cr}
$$
Note that
$C_1(x,y)=C_2(x,y)+C_3(x,y),\; C_4(x,y)=C_5(x,y)+C_6(x,y),\;
C_2(x,y)=\sign(C_5(x,y)),\; C_3(x,y)=\sign(C_6(x,y)).$

In the following list of objective functions (and of corresponding
methods) structure functions $C_1\cdc C_6$ are used not only for
quadratic assignment problems.

\medskip
\flushpar1. Three distinct extensions of the Slater (1961) method,
which had been originally suggested by Kendall (1955):  $\QA(R,C_k)$,
$k\in\{1,2,3\}$ (see, e.g., Arditti (1983)).

\medskip
\flushpar2. Three distinct extensions of the Kemeny (1959) method
(which is equivalent to the Slater method in the complete dichotomous
case):
\newline {\em minimize} $\;\;\suml_{(i,j,p)|\R}|r_{ij}^p-C_k(\r
i,\r j)|\;\;{\em for}\;\;\r\in\P$\newline
with $k\in\{1,2,3\}$ and summation over those $i,j,p$
for which $r_{ij}^p$ is defined in $\R$. According to Young (1986),
Kemeny's method had been initially proposed in a vague form by
Condorcet.

\medskip
\flushpar3. ``Weighted sum of back scores": $\QA(R,C_6).$ This method
was suggested by Thompson (1975) and Hubert (1976), and studied in Kano
and Sakamoto (1985), and Frey and Yehia-Alcoutlabi (1986).

\medskip
\flushpar4. ``Weighted sum of right scores": $\QA(R,C_5)$ (Kano and
Sakamoto (1983)).

\medskip
\flushpar5. ``Weighted sum of all scores": $\QA(R,C_4)$ (Chebotarev
(1988, 1990), Crow (1990)). This method can be reduced to ordering
alternatives by ``sum of wins minus sum of losses" (it is the Copeland
index; see Theorem~1 below) and is connected to some ideas of
Kendall (1970).

\medskip
\flushpar6. ``Net sum of back scores" (see, e.g., Weiss and Assous
(1987), Crow (1990)) $\QA((R-R^T)^{\smp},C_3)$:\newline
{\em maximize} $\suml_{i=1}^n\suml_{j:\: m^{}_{ij}>0}(r_{ij}-r_{ji})^{\smp}
C_3(\r i,\r j)$ {\em for} $\r\in\P.$

\medskip
\flushpar7. The following four methods are based on the idea of
balancing ``back scores" of two types: ``wins above" and ``losses
below" (Crow (1990, 1993)).

\medskip
\flushpar7a. Sum of absolute differences between Wins Above and Losses
Below -- ``WALB":
{\em minimize}
$\suml_{i=1}^n\Bigl\vert \suml_{j=1}^n(r_{ij}C_3(\r i,\r j)-
r_{ji}C_3(\r j,\r i))\Bigr\vert$
{\em for}
$\r\in\P.$

\flushpar7b. ``Refined WALB":
\newline
{\em minimize}
$\suml_{i=1}^n\Bigl\vert \suml_{j:\: m^{}_{ij}>0}{1\over m^{}_{ij}}
(r_{ij}C_3(\r i,\r j)- r_{ji}C_3(\r j,\r i))\Bigr\vert$
{\em for}
$\r\in\P.$

\medskip
\flushpar7c. ``Net WALB":\newline
{\em minimize}
$\suml_{i=1}^n\Bigl\vert \suml_{j:\: m^{}_{ij}>0}((r_{ij}-r_{ji})^{\smp}
C_3(\r i,\r j)-(r_{ji}-r_{ij})^{\smp}C_3(\r j,\r i))\Bigr\vert$
{\em for}
$\r\in\P.$

\medskip
\flushpar7d. ``Refined Net WALB":\newline
{\em minimize}
$\suml_{i=1}^n\Bigl\vert \suml_{j:\: m^{}_{ij}>0} {1\over m^{}_{ij}}
((r_{ij}-r_{ji})^{\smp}C_3(\r i,\r j)- (r_{ji}-r_{ij})^{\smp}C_3(\r
j,\r i))\Bigr\vert$
{\em for}
$\r\in\P.$

\medskip
\flushpar7e. ``Net-Difference-WALB":
\newline
{\em minimize}
$\suml_{i=1}^n\Bigl\vert \suml_{j:\: m^{}_{ij}>0}((r_{ij}-r_{ji})^{\smp}
C_6(\r i,\r j)- (r_{ji}-r_{ij})^{\smp}C_6(\r j,\r i))\Bigr\vert$
{\em for}
$\r\in\P.$

\medskip
The following method is new.

\medskip
\flushpar8. ``$\beta$-Least-Squares" ($\beta$-LS): {\em minimize}
$\suml_{(i,j,p)|\R}(r_{ij}^p- \beta C_4(\r i,\r j))^2$
{\em for}
$\r\in\P$\newline
with summation over those $i,j,p$ for which
$r_{ij}^p$ is defined in $\R$. Here $\beta$ is a positive real
parameter.

\medskip
{\it Remark.} In the methods based on ``net scores" (i.e., ``Net sum
of back scores", ``Net WALB", ``Refined Net WALB", and
``Net-Difference-WALB") a ``net draw" ($m\_{ij}>0,\;r\_{ij}=r\_{ji}$)
between two alternatives with different positions in $\r$ ($\r i\ne\r j$) is
worth being distinguished from the lack of comparisons between them
($m\_{ij}=0,$ where $r\_{ij}=r\_{ji}$ by definition). If only maximal
wins/losses are allowed, then the following modification provides
this distinction: replace $(r\_{ij}-r\_{ji})^{\smp}$ by
$\psi(r\_{ij}-r\_{ji})$, where

$$
\psi(z)=\cases{z,                     & if $z>0$,\cr
               (r_{\max}-r_{\min})/2, & if $z=0$,\cr
               0,                     & if $z<0$,\cr}
$$
as Crow (1990, 1993) proposes for the case $r_{\max}=1,\:$
$r_{\min}=0.$ Note that such a modification preserves our results that
involve these methods, i.e., Corollary~1 and Theorem~6 below. Another
possible modification based on the function
$\psi'(z)=(z+r_{\max}-r_{\min})^{\smp}$ preserves Corollary~1 and
Theorem~6 as well. Above we wrote $j:\:  m^{}_{ij}>0$ under the second
sums in the objective functions of the ``Net"-methods in order to
support these possible modifications.
\medskip

Some other methods can be obtained by extending the measures of
association from Critchlow (1985) to incomplete paired comparisons.

\section{Generalized row sum method}

The generalized row sum method (Chebotarev (1989, 1994)) is not based
on the resolution of a discrete optimization problem, however it has
some connection with the $\beta$-LS method (Theorem~3 below). On the
other hand, the generalized row sum method will be shown to satisfy
Self-Consistent Monotonicity, the main axiom in this paper
(Theorem~8).

For the sake of simplicity, we suppose here that incomplete paired
comparison matrices $R^{(p)},$ $p=1\cdc m$, are skew-symmetric: if
$r_{ij}^p$ is defined in $R^{(p)}$, then $r_{ji}^p=-r_{ij}^p.$ In this
case, $r_{\min}=-r_{\max}$. The generalized row sum method estimates
the alternatives by the indexes $x_1\cdc x_n$ ({\it generalized row sums})
that satisfy the following system of linear equations:
$$
x\_i=\suml_{(k,p)|i}\big(r_{ik}^p+\ve\cdot(x\_k-x\_i+r_{ik}^pmn)\big),
\quad i=1\cdc n,
\eqno(4)
$$
where $\ve$ is a nonnegative parameter. This system of equations
has been proven to have a unique solution for every $\R.$  The
corresponding optimal weak order $\r$ is defined as follows:
$X_i\gro X_j$ iff $x\_i>x\_j.$

The generalized row sum method is an extension of the row sum method
(and of the Borda rule in the case where individual preferences are
linear orders) to incomplete paired comparisons. Specifically, if
$\R$ is complete, then for any $\ve\ge0$, $\:x\_i=s\_i\;\:(i=1\cdc
n)$ holds, where
$$
s\_i=\suml_{(k,p)|i}r_{ik}^p=t\_i/2.
$$

This method has been derived both axiomatically and statistically. The
value $f_{ik}^p=r_{ik}^p+\ve\cdot(x\_k-x\_i+r_{ik}^pmn)$ is the
contribution of the comparison outcome $r_{ik}^p$ to the estimate
$x\_i$ of $X_i$. Parameter $\ve\ge0$ is said to be {\it reasonable} for
given $n$ and $m$ if for any array $\R$ that consists of $m$ $n$-by-$n$
paired comparison matrices, the value
$$
f_{ik}^p=r_{ik}^p+\ve\cdot(x\_k-x\_i+r_{ik}^pmn)
$$
is non-negative at $r_{ik}^p=r_{\max}$ (maximal win) and
non-positive at $r_{ik}^p=r_{\min}=-r_{\max}$ (maximal loss), for
any $i,$ $j,$ and $p.$

It has been shown that the reasonableness of $\ve$ is equivalent
to satisfying the constraint
$$
0\le\ve\le{1\over m(n-2)}.
$$

\section{Some connections to direct methods}

In this section, we prove three theorems concerning connections between
discrete optimization methods, namely $\QA(R,C_4)$ and $\beta$-LS, and
direct methods for aggregating preferences. The first two theorems are
formulated for the general case of weak orders ($\P=\W$). Theorem~1
shows that the method $\QA(R,C_4)$ can be reduced to ordering alternatives
in the decreasing order of their Copeland indexes (with an arbitrary
order of the alternatives having the same Copeland index). Note that the
related problems $\QA(R,C_5)$ and $\QA(R,C_6)$ are, in general,
NP-complete.

\bigskip\noindent
{\sc Theorem 1} (Reduction of $\WQA(R,C_4)$ to the Copeland ranking):
{\it A weak order $\r$ is a solution of $\WQA(R,C_4)$ for $\R$ if and
only if}
\vspace{-4pt}
$$
{\em for\; any\ } X_i,X_j\in\X,\quad t_i>t_j {\em\ implies\ }
X_i\gro X_j.
\eqno(5)
$$

The proofs of all statements are given in the Appendix. An analogous
theorem for {\it linear} orders has been proved in Chebotarev (1988,
1990).

A similar statement holds for the $\beta$-LS method with small enough
$\beta$.

\bigskip\noindent
{\sc Theorem 2} (Partial reduction of $\beta$-LS with small $\beta$ to
the Copeland ranking): {\it Let $\R$ be an array of paired comparisons
on $\X$. There exists a number $\beta_0>0$ such that: if
$\;0<\beta<\beta_0$, then every solution $\r^*\in\W$ of the $\beta$-LS
problem with parameter $\beta$ for $\R$ satisfies the following
condition:}
\vspace{-4pt}
$$
{\em for\; any\ } X_i,X_j\in\X,\quad t_i>t_j {\em\ implies\ }
X_i\succ_{\r^*} X_j.
\eqno(6)
$$

There is an important difference between the methods $\QA(R,C_4)$ and
$\beta$-LS with a small $\beta$. Namely, according to Theorem~2, (6) is
a necessary but not a sufficient condition of optimality. In other
words, the latter method does not permit arbitrariness in ordering the
alternatives with the same Copeland index. Indeed, it can be easily shown by
examples that $\beta$-LS with a small parameter may yield a narrower
set of optimal orders than $\QA(R,C_4)$.

Now consider a continuous counterpart of the $\beta$-Least-Squares
method. Note that for any linear order $\r$ on $\X$,
$\suml_{i=1}^n\r i=0$ and $\suml_{i=1}^n(\r i)^2= {1\over
3}(n-1)n(n+1)$. Denote the latter value by $D_n^2$ and consider the
following {\it relaxed $\beta$-LS method:}
$$
{\em minimize} \;\;\suml_{(i,j,p)|\R}(r_{ij}^p-
\beta(y_i-y_j))^2\;\;{\em for\;real}\;\;y_1\cdc y_n
\eqno(7)
$$
subject to
\vspace{-4pt}
$$
\suml_{i=1}^n y_i=0
\eqno(8)
$$
and
\vspace{-4pt}
$$
\suml_{i=1}^n y_i^2=D_n^2.
\eqno(9)
$$

The difference between $\beta$-LS and relaxed $\beta$-LS is that for
the former problem the set of admissible solutions is narrower: not the
whole intersection of the hyperplane $\suml_{i=1}^n y_i=0$ with the
hypersphere $\suml_{i=1}^n y_i^2=D_n^2$, but the set of points obtained
from $(-(n-1),-(n-3)\cdc(n-1))$ by all possible permutations of the
coordinates. These points are all vertices of a specific polyhedron
(polytope) inscribed into that intersection. According to the
following theorem, relaxed $\beta$-LS is closely connected to the
generalized row sum method.

\bigskip\noindent
{\sc Theorem 3} (Reduction of the relaxed $\beta$-LS to the generalized
row sums): {\it Let $y=(y_1\cdc y_n)$ be a solution of the relaxed
$\beta$-LS problem with some $\beta$ for an array of paired comparisons
$\R=(r_{ij}^p)_{i,j\in\{1\cdc n\}}^{p\in\{1\cdc m\}}$. Let $\R'$ be an
array of paired comparisons with elements $(r_{ij}^p)'
=r_{ij}^p-r_{ji}^p$. Then for some $\ve$, the vector $y$ is
proportional to the vector of generalized row sums obtained with
parameter $\ve$ for $\R'$.}
\bigskip

Another example of Lagrangian relaxation applied to a discrete
preference aggregation problem can be found in Arditti (1983).

\section{Self-Consistency and Self-Consistent\newline Monotonicity}

H.A.~David (1987) said ``...nonparametric method cannot be entirely
satisfactory when the $m_{ij}$ differ greatly." Our aim is to
investigate to what extent such a method can be satisfactory, and so
we examine the properties of the methods above. In this section, a
new axiom named {\it Self-Consistency\/} and its extension, {\it
Self-Consistent Monotonicity} are introduced.

Let us say that an outcome $(r_{ik}^p, r_{ki}^p)$ of comparing $X_i$ to
$X_k$ is {\it not weaker with respect to a weak order\/} $\r$ than an
outcome $(r_{j\ell}^q, r_{\ell j}^q)$ of comparing $X_j$ to $X_\ell$
iff $r_{ik}^p\ge r_{j\ell}^q$, $r_{ki}^p\le r_{\ell j}^q$, and
$X_k\gero X_\ell$. If, in addition, at least one of the inequalities
(relations) is strict, then the outcome $(r_{ik}^p, r_{ki}^p)$ is said
to be {\it stronger\/} than $(r_{j\ell}^q, r_{\ell j}^q)$ with respect
to $\r$.

\medskip
{\bf Self-Consistency.} For any optimal weak order $\r$ and for any
$X_i, X_j\in\X$, the statement [There exists a one-to-one
correspondence between the set of comparison outcomes of $X_i$ and the
set of comparison outcomes of $X_j$ such that each outcome of $X_i$ is
not weaker than the corresponding outcome of $X_j$ with respect to
$\r$] implies [$X_i\gero X_j$]. If, in addition, at least one outcome
of $X_i$ is stronger than the corresponding outcome of $X_j$ with
respect to $\r$, then $X_i\gro X_j$.
\medskip

Self-Consistency enables us to confront two alternatives having the same
number of comparisons. Now suppose that the alternative dominating in such a
confrontation achieves several extra {\it maximal wins\/} and the
dominated alternative gets some number of extra {\it maximal losses}. It
is reasonable to demand that this addition of extra outcomes
preserves the result of confrontation: the former alternative remains
``better". Let us extend Self-Consistency in this way.  \medskip

{\bf Self-Consistent Monotonicity (SCM).} Suppose $\r$ is an optimal
weak order and $X_i,X_j\in\X$. Let $\R_i$ and $\R_j$ be the sets of
comparison outcomes of $X_i$ and $X_j$, respectively. Suppose that
$\R_i=\R_i'\cup\R_i''$\ $(\R_i\cap\R_i''=\emptyset)$,
$\R_j=\R_j'\cup\R_j''$\ $(\R_j\cap\R_j''=\emptyset)$, $\R_i''$
consists of maximal wins, $\R_j''$ consists of maximal losses, and
there exists a one-to-one correspondence between $\R_i'$ and $\R_j'$
(in particular, $\R_i'$ and $\R_j'$ may be empty) such that every
outcome from $\R_i'$ is not weaker than the corresponding outcome from
$\R_j'$ with respect to $\r$. Then $X_i\gero X_j$. If, in
addition, at least one outcome from $\R_i'$ is stronger than the
corresponding outcome from $\R_j'$ with respect to $\r$ or
$\R_i''\ne\emptyset$ or $\R_j''\ne\emptyset$, then $X_i\gro X_j$.
\medskip

Possibly, some analysts can be inclined to consider the entire set of
optimal orders as an indivisible macro-decision whose elements
represent different characteristic features of the set of original
preferences. From this point of view, optimal orders should be
considered not separately but jointly, and Self-Consistency which
addresses to every separate optimal order is a surplus requirement. A
possible objection to this opinion is as follows. In most situations
we have to make only one decision.  As soon as it is made, any
appealing to other optimal decisions becomes out of place. The
decision we make should be logical by itself, apart from rejected
opportunities.

\section{All strict operators break Self-Consistency}
\smallskip

Recall that an aggregating operator is {\it strict\/} if its optimal
orders are always linear.
\smallskip

\noindent
{\sc Theorem 4}: {\it If an aggregating operator is strict,
then it does not satisfy Self-Consistency.}
\bigskip

This theorem has an easy but somewhat degenerate proof. Indeed, note
that Self-Consistency does not prohibit the sets of comparison outcomes
of $X_i$ and $X_j$ to be empty. In this case, Self-Consistency implies
$X_i\gero X_j$ and $X_j\gero X_i$, which is broken by any linear order.
A similar proof with alternatives that have nonempty sets of comparisons and
$n>2$ can be carried out by considering the following $\R$:
$r_{13}^1=r_{23}^1=r_{\max},\; r_{31}^1=r_{32}^1=r_{\min};$ all other
$r_{ij}^p$ with $i\ne j$ are undefined. In the Appendix we give another
proof, which demonstrates the application of Self-Consistency to cyclic
preferences.
\bigskip

\section{Operators generating weak orders}

Theorem 4 motivates the consideration of aggregating operators that
generate not only linear orders but arbitrary weak orders. In
particular, we shall consider {\it weak quadratic assignment
problems\/} $\WQA(R,C)$, i.e., problems (3) with $\P=\W$ (note that
Theorem~1 and Theorem~2 have been formulated for this general case).

To that end it is useful to modify structure functions $C_2$, $C_3$,
$C_4$, and $C_6$.  Indeed, note that the structure functions
$C_1(x,y)\cdc C_6(x,y)$ depend on $x-y$. Suppose $g(d)$ is the
contribution of the comparison outcome $(r_{ij}^p=r_{\max},\;
r_{ji}^p=r_{\min})$ to the quadratic assignment objective function,
provided that $\r i-\r j=d$. Then, by (3), $g(d)=
r_{\max}C(d,0)+r_{\min}C(0,d)$. It is reasonable to require
$$
g(-1)<g(0)<g(1).
\eqno(10)
$$

Indeed, since the quadratic assignment objective function measures
consistency between the original preferences and a tentative
resulting order, this requirement is motivated by that maximal win is
more natural for alternatives with higher social estimate.

For $C_1$ and $C_4$, (10) amounts to $r_{\max}>r_{\min}$,
whereas for $C_2,C_3,C_5$, and $C_6$ it is equivalent to [$r_{\max}>0$
and $r_{\min}<0$]. Therefore (10) is broken even for the customary
sporting point systems: $r_{ij}^p\in\{0,{1\over2},1\}$ and
$r_{ij}^p\in\{0,1,2\}$. As a result, for these point systems, the
weak order in which all alternatives are tied is never
optimal for $\QA(R,C_2)$ and $\QA(R,C_5)$, and is always optimal for
$\QA(R,C_3)$ and $\QA(R,C_6)$.

Thus let us revise $C_2,C_3,C_5$, and $C_6$ as follows:
$$
C'_2(x,y)=\sign(x-y)+1,\quad C'_3(x,y)=\sign(x-y)-1,
$$
$$
C'_5(x,y)=(x-y+1)^{\smp},\quad C'_6(x,y)=(x-y-1)^{\smm}.
$$

$C_1$ and $C_4$ do not require revisions: let
$C'_1(x,y)=C_1(x,y),\;\:C'_4(x,y)=C_4(x,y).$

For all these functions, (10) amounts to $r_{\max}>r_{\min}$, and
they are equivalent to their prototypes in all optimization
methods of Section~4 in the strict case. In the rest of the paper, we
consider quadratic assignment problems and other problems of
Section~4 with $C'_k$ substituted for $C_k$ $(k\in\{1\cdc 6\})$ and
$\P=\W$.

\section{Indifference to the degree of resulting\newline
         preferences contradicts SCM}

Let us say that an aggregating operator {\it equalizes\/} weak
orders $\r$ and $\r'$ for $\R$ if $\r$ and $\r'$ are both optimal
for $\R$ or both are not optimal. An aggregating operator will be
called {\it indifferent to the degree of resulting preferences\/}
if it equalizes every $\r$ and $\r'$ such that
$$
{\em for\ all\ }\; r_{ij}^p \;{\em defined\ in\ } \R,\;\;\sign(\r i-\r
j) =\sign(\r'i-\r'j).
$$

\noindent
{\sc Theorem 5}: {\it If a nonstrict aggregating operator is
indifferent to the degree of resulting preferences and $n>2$, then it
violates { SCM.}}
\bigskip

{\sc Corollary 1}: {\it The nonstrict aggregating operators
corresponding to:\newline $\WQA(R,C'_k),$ $k\in\{1,2,3\};$
{\em minimize}\hspace{-3mm} $\suml_{(i,j,p)|\R}|r_{ij}^p-C'_k(\r i,\r
j)|\;\;{\em for}\;\r\in\W,$ $k\in\{1,2,3\}$ (extensions of the
Kemeny median); ``Net sum of back scores", ``WALB", ``Net WALB",
``Refined WALB", and ``Refined Net WALB" violate SCM.}

\section{Are there discrete optimization methods\newline that obey
Self-Consistent Monotonicity?}

\noindent
{\sc Theorem 6:} {\it If $n>2$, then the nonstrict aggregating
operators corresponding to $\WQA(R,C'_k)$ with $k\in\{4,5,6\}$ and
``Net-Difference-WALB" violate SCM.}
\bigskip

The claim that the $\beta$-LS operator satisfies SCM might provide
a ``happy end" of this paper. However, this is not the case.

\bigskip\noindent
{\sc Theorem 7:} {\it If $n>4$ then the $\beta$-LS operator violates
Self-Consistency for any $\beta>0$.}
\bigskip

Recall that the $\beta$-LS method can be considered as a discrete
analog of the generalized row sum method (Theorem~3).

\bigskip\noindent
{\sc Theorem 8:} {\it The generalized row sum method with positive
$\ve$ satisfies Self-Consistency. Moreover, it satisfies
Self-Consistent Monotonicity when $\ve$ is positive and
reasonable.}
\bigskip

Comparison of Theorem~7 and Theorem~8 suggests that the $\beta$-LS
method fails to satisfy Self-Consistency because of its discreteness.
Indeed, in the proof of Theorem~7 given in the Appendix, $X_1$ has a
small superiority over $X_2$ in the original preferences, and
Self-Consistency requires $X_1\gro X_2.$ However, $\beta$-LS ties $X_1$
and $X_2$ for every $n>5.$ Note that $\beta$-LS minimizes some kind of
proximity between the initial preferences and the tested weak orders.
The superiority of $X_1$ over $X_2$ turns out to be so small that it is
closer to ``draw" than to ``win". This is typical of nonstrict
discrete methods like $\beta$-LS. There are only three possible
relations between two alternatives in a social weak order, ``worse",
``better", and ``equivalent", and the latter turns out to be optimal
for small superiorities under the nonstrict aggregating procedures.
Is this a shortcoming or not?  We believe that in case we must choose
only one alternative, even a small superiority is worth being taken
into account, and so such a tie is not useful.  An advantage of
continuous approaches is that they enable one to measure the
differences between the adjacent alternatives, whereas the discrete
methods give no means for that. (However, some information can be
extracted through comparing the optimal value of the objective
function with its values for orders where these alternatives are tied
or interchanged.)

All the discrete optimization methods we considered proved to break
Self-Consistent Monotonicity. Nevertheless, the question in the heading
of this section is a methodological rather than a mathematical one.
Indeed, a discrete optimization method that satisfies SCM can be
designed artificially, for example, by using explicit expressions
of the generalized row sums $x\_1\cdc x\_n$:
$$
{\em maximize}
\;\;\suml_{i=1}^n (x\_i\r i-\alpha(\r i)^2) \;\;\;{\em for}\;\r\in\W,
\eqno(11)
$$
where $\alpha$ is a small enough positive constant.

To prove that the optimal values $\r i$ are ordered exactly as $x\_i$,
note that every maximizing weak order for the objective function
$\suml_{i=1}^n x\_i\r i$ preserves the strict component of the order of
$x\_1,...,x\_n$ (Lemma~1 in the proof of Theorem~1), and the
subtraction of $\alpha(\r i)^2$ in (11) provides equal $\r i$ for
the alternatives with equal generalized row sums $x\_i$. (Indeed,
equal numbers provide a minimum for the sum of squares subject to
their fixed sum.) It follows that aggregating operator (11)
satisfies SCM.  However, such a method would remain essentially
based on ``continuous" indexes.  Now we do not know any proper
discrete optimization operators that satisfy SCM.

\section{Conclusion}

If an aggregating operator is strict, then it breaks Self-Consistent
Monotonicity (SCM), since this axiom prescribes equivalence of some
alternatives (Theorem~4). Many aggregating operators associated with
discrete optimization problems are ``indifferent to the degree of
resulting preferences", which is incompatible with SCM (Theorem~5).
Nonstrict discrete optimization methods like $\beta$-LS violate
Self-Consistent Monotonicity, since they produce equivalence of some
alternatives, one of which having a small superiority over another. On the
other hand, there are ``continuous" methods that satisfy SCM, for
example, the generalized row sum method (Theorem~8).

The transfer from linear orders to weak orders is the first step of
relaxation. Possibly, this step is not sufficient for such a keen type
of data as unbalanced (incomplete) preferences. A next possible step is
the conversion to aggregation models with real unknown parameters
that measure the value (utility) of alternatives. Such {\it indirect
scoring procedures} are considered in Chebotarev and Shamis
(1996) where a sufficient condition of Self-Consistent
Monotonicity and some more positive examples are given.

\section*{Appendix: Proofs}

{\sc Proof of Theorem 1:} Suppose $\r$ is an arbitrary weak order on
$\X$, and $f(\r)$ is the value of the objective function for $\r$. Then
$$ f(\r)= \suml_{i=1}^n\suml_{j=1}^n r_{ij}(\r i-\r j)=
\suml_{i=1}^n\left(\r i \suml_{j=1}^n r_{ij}\right)-
\suml_{j=1}^n\left(\r j \suml_{i=1}^n r_{ij}\right)=
\suml_{i=1}^n\r i\,t_i.
\eqno(A1)
$$

Now it suffices to prove the following lemma.

\medskip\noindent
{\sc Lemma 1} (A weak order maximizes scalar product iff it preserves
relation ``$>$"): {\it For any real vector $u=(u\_1\cdc u\_n),$ a weak
order $\r$ is a solution of the problem}
\vspace{-4pt}
$$
{\em maximize}
\;\;\suml_{i=1}^n u\_i\r i \;\;\;{\em for}\;\r\in\W
\eqno(A2)
$$
{\em if and only if}
\vspace{-4pt}
$$
{\em for\; any\ } X_i,X_j\in\X,\quad u_i>u_j {\em\ implies\ }
X_i\gro X_j.
\eqno(A3)
$$

\medskip
{\sc Proof of Lemma 1}: Let $\r$ be a solution of the problem (A2).
Assume that there exist $X_k$ and $X_\ell$ such that $u\_i>u\_j$, but $\r
k\le \r \ell$. Consider two cases.

(A) $\r k< \r \ell$. Consider the weak order $\r'$ that is obtained
from $\r$ by interchanging $X_k$ and $X_\ell$. Then, according to (A1),
$$
f(\r)-f(\r')
=\suml_{i=1}^n(\r i-\r'i)u\_i
=(\r k)u\_k+(\r \ell)u\_\ell-(\r'k)u\_k-(\r'\ell)u\_\ell
$$
$$
=(\r k)u\_k+(\r \ell)u\_\ell-(\r \ell)u\_k-(\r k)u\_\ell
=(\r k-\r \ell)(u\_k-u\_\ell)<0,
$$
and $\r$ cannot be a solution of $\QA(R,C_4)$, in contradiction to
the assumption.

(B) $\r k=\r\ell$. Let $\X_{k\ell}=\{X_j:\r j=\r k\}\setminus
\{X_k,X_\ell\}$. Recall that $\r$ is a binary relation, i.e., $X_v\gero
X_w$ is a designation of $(X_v,X_w)\in\r$. We shall ``move apart"
$X_k$ and $X_\ell$ in $\r$ preserving the positions of all other
alternatives. Consider the weak order $\r'$ that is obtained by removing
from $\r$ the pair $(X_\ell,X_k)$ and the pairs $(X_j,X_k)$ and
$(X_\ell,X_j)$ for all $X_j\in\X_{k\ell}$. Then $\r' k=\r k
+|\X_{k\ell}|+1,\:$ $\r' \ell=\r \ell-|\X_{k\ell}|-1$, and $\r i'=\r i$
for all $X_i\in\X\setminus\{X_k,X_\ell\}$. Hence, by (A1),
$$
f(\r)-f(\r')
=\suml_{i=1}^n(\r i-\r'i)u\_i
=(-|\X_{k\ell}|-1)u\_k+(|\X_{k\ell}|+1)u\_\ell
$$
$$
=(u\_\ell-u\_k)(|\X_{k\ell}|+1)<0,
$$
and $\r$ cannot be a solution of $\QA(R,C_4)$, in contradiction to
our assumption. Necessity of (A3) is shown.

To prove sufficiency of (A3), note that all weak orders satisfying (A3)
can be obtained one from another by sequential adding and removing
pairs $(X_k,X_\ell)$ such that $u\_k=u\_\ell$. Therefore, if for an
arbitrary $i$ we denote the set $\{X_j\in\X:u\_j=u\_i\}$ by $\X_i$,
then the value $\suml_{X_j\in\X_i}\r j $ is the same for all weak
orders satisfying (A3). Consequently, all such weak orders have the
same (and thus maximal!) value of (A1). This completes the proof of
Lemma~1 and Theorem~1.
\bigskip

{\sc Proof of Theorem 2:} Suppose $\r^*$ satisfies (6) and let for any
$\r$
\vspace{-4pt}
$$
f(\r)=
\suml_{(i,j,p)|\R}(r_{ij}^p- \beta C_4(\r i,\r j))^2
$$
$$
=\suml_{(i,j,p)|\R}(r_{ij}^p)^2
-2\beta\suml_{i=1}^n\suml_{j=1}^n r_{ij}C_4(\r i, \r j)
+\beta^2\hspace{-1.7ex}\suml_{(i,j,p)|\R}(\r i- \r j)^2.
\eqno(A4)
$$
The second term in the right-hand side of (A4) corresponds to
$\QA(R,C_4)$. By Theorem~1, there exists a number $E>0$ such that for
any $\r'$ not satisfying (6),
$$
 \suml_{i=1}^n\suml_{j=1}^n r_{ij}C_4(\r^* i, \r^* j)
-\suml_{i=1}^n\suml_{j=1}^n r_{ij}C_4(\r'  i, \r'  j)
>E.
$$
Let $F$ be an upper bound of $\suml_{(i,j,p)|\R}(\r i- \r j)^2$ for all
weak orders $\r$ on $\X$. Then for any $\beta$ such that
$0<\beta<{2E\over F}=\beta_0$ and for any $\r'$ not satisfying (6),
$f(\r^*)>f(\r')$, and $\r'$ cannot be a solution. The theorem is
proved.
\bigskip

{\sc Proof of Theorem 3:} This proof is technical, and we give only a
plan. Searching the minimum (7) subject to (8) and (9) with the
Lagrange multiplier method, we get a system of $n$ linear equations
in $y_1\cdc y_n$. Summing all of them, we derive that the multiplier
corresponding to (8) equals zero and then conclude that this system of
equations for some $\ve$ coincides with that of the generalized
row sum method. This completes the proof.
\bigskip

{\sc Proof of Theorem 4:} This proof is very simple, however we give it
in detail in order to illustrate the application of Self-Consistency.
Assume on the contrary that there exists a strict aggregating operator
that satisfies Self-Consistency.

 1. Let $n>2$. Consider the following array of paired comparisons $\R$:
$r_{12}^1=r_{23}^1=r_{31}^1=r_{\max};\; r_{21}^1 =r_{32}^1 =r_{13}^1
=r_{\min};$ all other $r_{ij}^p$ with $i\ne j$ are undefined (Figure
A1).

\begin{figure}[htb]
\setlength{\unitlength}{5mm}
\begin{picture}(26,13.6)
\put( 2, 5){\line(1,0){11}}
\put( 2,12){\line(1,0){11}}
\put( 2, 5){\line(0,1){7}}
\put(13, 5){\line(0,1){7}}
\put(2,5){\begin{picture}(13,7)
             \put(1.5,5.0){\circle{0.2}} 
             \put(5.5,4.0){\circle{0.2}} 
             \put(2.5,2.0){\circle{0.2}} 
             \put(7.0,3.0){\circle{0.2}} 
             \put(9.5,3.0){\circle{0.2}} 
             \put(1.2,5.3){$X_1$}
             \put(5.2,4.3){$X_2$}
             \put(2.1,0.9){$X_3$}
             \put(6.6,1.9){$X_4$}
             \put(9.1,1.9){$X_n$}
             \put(7.9,1.9){$\cdots$}
             \put(1.8,4.89) {\vector(4,-1) {3.2}}  
             \put(5.2,3.8)  {\vector(-3,-2){2.4}}  
             \put(2.35,2.3) {\vector(-1,3) {0.75}} 
             \put(5.0,-1){$R^{(1)}$}
          \end{picture}}
\put(16, 5){\line(1,0){11}}
\put(16,12){\line(1,0){11}}
\put(16, 5){\line(0,1){7}}
\put(27, 5){\line(0,1){7}}
\put(16,5){\begin{picture}(13,7)
             \put(1.5,5.0){\circle{0.2}} 
             \put(5.5,4.0){\circle{0.2}} 
             \put(2.5,2.0){\circle{0.2}} 
             \put(7.0,3.0){\circle{0.2}} 
             \put(9.5,3.0){\circle{0.2}} 
             \put(1.2,5.3){$X_1$}
             \put(5.2,4.3){$X_2$}
             \put(2.1,0.9){$X_3$}
             \put(6.6,1.9){$X_4$}
             \put(9.1,1.9){$X_n$}
             \put(7.9,1.9){$\cdots$}
             \put(2.5,-1){$R^{(2)}=\cdots=R^{(m)}$}
          \end{picture}}
\put(2.0,2.0){\small{\sc Figure A1.}  $\R=(R^{(1)}\cdc R^{(m)})$
in the proof of Theorem~4: $n>2.$}
\put(2.0,1.0){\small Every arc $(X_i,X_j)$ designates that
$r_{ij}^p=r_{\max}$ and $r_{ji}^p=r_{\min}.$}
\end{picture}
\end{figure}

Let $\r$ be a linear order on $\X$ in which $X_1\gro X_2$ and $X_2\gro
X_3$.  The set of comparison outcomes of $X_1$ is
$\{(r_{12}^1=r_{\max},\:  r_{21}^1=r_{\min}),\: (r_{13}^1=r_{\min},\:
r_{31}^1=r_{\max})\};$ the set of comparison outcomes of $X_3$ is
$\{(r_{31}^1=r_{\max},\:  r_{13}^1=r_{\min}),\: (r_{32}^1=r_{\min},\:
r_{23}^1=r_{\max})\}$.  Note that $(r_{31}^1, r_{13}^1)$ is stronger
than $(r_{12}^1, r_{21}^1)$ with respect to $\r$, since $X_1\gro X_2$,
and $(r_{32}^1, r_{23}^1)$ is stronger than $(r_{13}^1, r_{31}^1)$ with
respect to $\r$, since $X_2\gro X_3$. Therefore $X_1\gro X_3$
contradicts Self-Consistency, and $\r$ is not optimal. Arguing as above
we obtain that any linear order for which $X_2\gro X_3\gro X_1$ or
$X_3\gro X_1\gro X_2$ is not optimal. Now assume $X_1\gro X_3$ and
$X_3\gro X_2$. Let us compare the outcomes of $X_2$ and $X_1$.  Note
that $(r_{23}^1, r_{32}^1)$ is stronger than $(r_{12}^1, r_{12}^1)$
with respect to $\r$, since $X_3\gro X_2$ and $(r_{21}^1, r_{12}^1)$ is
stronger than $(r_{13}^1, r_{31}^1)$, since $X_1\gro X_3$. Therefore,
$X_1\gro X_2$ contradicts Self-Consistency, and $\r$ is not optimal. In
the same way we conclude that if $X_2\gro X_1\gro X_3$ or $X_3\gro
X_2\gro X_1$, then $\r$ is not optimal. Since every linear order on
$\X$ obeys one of the above six assumptions, we obtain that the set of
optimal orders is empty in contradiction to the definition of
aggregating operator.

2. Let $m>1$. Consider the following $\R$: $r_{12}^1=r_{21}^2=r_{\max};\;
r_{21}^1=r_{12}^2=r_{\min};$ all other $r_{ij}^p$ with $i\ne j$ are
undefined (Figure A2).

\begin{figure}[htb]
\setlength{\unitlength}{5mm}
\begin{picture}(26,12.6)
\put(2,5){\line(1,0){7}}
\put(2,11){\line(1,0){7}}
\put(2,5){\line(0,1){6}}
\put(9,5){\line(0,1){6}}
\put(2,5){\begin{picture}(7,6)
             \put(1.5,4.0){\circle{0.2}} 
             \put(5.5,4.0){\circle{0.2}} 
             \put(2.0,2.0){\circle{0.2}} 
             \put(5.0,2.0){\circle{0.2}} 
             \put(1.2,4.3){$X_1$}
             \put(5.2,4.3){$X_2$}
             \put(1.6,0.9){$X_3$}
             \put(4.6,0.9){$X_n$}
             \put(3.1,0.9){$\cdots$}
             \put(1.8,4.0){\vector(1,0){3.45}} 
             \put(3.0,-1){$R^{(1)}$}
          \end{picture}}
\put(11,5){\line(1,0){7}}
\put(11,11){\line(1,0){7}}
\put(11,5){\line(0,1){6}}
\put(18,5){\line(0,1){6}}
\put(11,5){\begin{picture}(7,6)
             \put(1.5,4.0){\circle{0.2}} 
             \put(5.5,4.0){\circle{0.2}} 
             \put(2.0,2.0){\circle{0.2}} 
             \put(5.0,2.0){\circle{0.2}} 
             \put(1.2,4.3){$X_1$}
             \put(5.2,4.3){$X_2$}
             \put(1.6,0.9){$X_3$}
             \put(4.6,0.9){$X_n$}
             \put(3.1,0.9){$\cdots$}
             \put(5.2,4.0){\vector(-1,0){3.45}} 
             \put(3.0,-1){$R^{(2)}$}
          \end{picture}}
\put(20,5){\line(1,0){7}}
\put(20,11){\line(1,0){7}}
\put(20,5){\line(0,1){6}}
\put(27,5){\line(0,1){6}}
\put(20,5){\begin{picture}(7,6)
             \put(1.5,4.0){\circle{0.2}} 
             \put(5.5,4.0){\circle{0.2}} 
             \put(2.0,2.0){\circle{0.2}} 
             \put(5.0,2.0){\circle{0.2}} 
             \put(1.2,4.3){$X_1$}
             \put(5.2,4.3){$X_2$}
             \put(1.6,0.9){$X_3$}
             \put(4.6,0.9){$X_n$}
             \put(3.1,0.9){$\cdots$}
             \put(0.55,-1){$R^{(3)}=\cdots=R^{(m)}$}
          \end{picture}}
\put(2.0,2.0){\small{\sc Figure A2.}  $\R=(R^{(1)}\cdc R^{(m)})$
in the proof of Theorem~4: $m>1.$}
\put(2.0,1.0){\small Every arc $(X_i,X_j)$ designates that
$r_{ij}^p=r_{\max}$ and $r_{ji}^p=r_{\min}.$}
\end{picture}
\end{figure}

In the same way as above, we obtain that any optimal linear
order can contain neither $(X_1,X_2)$ nor $(X_2,X_1)$, and the set of
optimal orders is empty. Again we have contradiction with the
definition of aggregating operator.

3. We have not covered the case $m=1,\: n=2$ yet. Here we can only
offer the degenerate proof described in Section~8. If {\it draws\/}
(i.e., $r_{ij}^p=r_{ji}^p$) are allowed, this provides a more
sensible proof.
\bigskip

{\sc Proof of Theorem 5:} Consider any aggregating operator that is
indifferent to the degree of resulting preferences. Assume that it
satisfies SCM.

1. Let $n>3$. Consider the following $\R:$
$r_{13}^1=r_{34}^1=r_{24}^1=r_{\max},$
$r_{31}^1=r_{43}^1=r_{42}^1=r_{\min};$ all other $r_{ij}^p$ with $i\ne
j$ are undefined (Figure A3).

\begin{figure}[htb]
\setlength{\unitlength}{5mm}
\begin{picture}(26,14.6)
\put( 2, 5){\line(1,0){11}}
\put( 2,13){\line(1,0){11}}
\put( 2, 5){\line(0,1){8}}
\put(13, 5){\line(0,1){8}}
\put(2,5){\begin{picture}(13,7)
             \put(1.5,6.0){\circle{0.2}} 
             \put(1.5,4.0){\circle{0.2}} 
             \put(5.5,5.0){\circle{0.2}} 
             \put(5.5,2.0){\circle{0.2}} 
             \put(7.2,4.0){\circle{0.2}} 
             \put(9.7,4.0){\circle{0.2}} 
             \put(1.2,6.3){$X_1$}
             \put(1.1,3.0){$X_3$}
             \put(5.2,5.3){$X_2$}
             \put(5.2,1.0){$X_4$}
             \put(6.8,3.0){$X_5$}
             \put(9.3,3.0){$X_n$}
             \put(8.1,3.0){$\cdots$}
             \put(1.5,5.7) {\vector(0,-1){1.4}} 
             \put(1.8,3.87){\vector(2,-1){3.3}} 
             \put(5.5,4.7) {\vector(0,-1){2.3}} 
             \put(5.0,-1){$R^{(1)}$}
          \end{picture}}
\put(16, 5){\line(1,0){11}}
\put(16,13){\line(1,0){11}}
\put(16, 5){\line(0,1){8}}
\put(27, 5){\line(0,1){8}}
\put(16,5){\begin{picture}(13,7)
             \put(1.5,6.0){\circle{0.2}} 
             \put(1.5,4.0){\circle{0.2}} 
             \put(5.5,5.0){\circle{0.2}} 
             \put(5.5,2.0){\circle{0.2}} 
             \put(7.2,4.0){\circle{0.2}} 
             \put(9.7,4.0){\circle{0.2}} 
             \put(1.2,6.3){$X_1$}
             \put(1.1,3.0){$X_3$}
             \put(5.2,5.3){$X_2$}
             \put(5.2,1.0){$X_4$}
             \put(6.8,3.0){$X_5$}
             \put(9.3,3.0){$X_n$}
             \put(8.1,3.0){$\cdots$}
             \put(2.5,-1){$R^{(2)}=\cdots=R^{(m)}$}
          \end{picture}}
\put(2.0,2.0){\small{\sc Figure A3.}  $\R=(R^{(1)}\cdc R^{(m)})$
in the proofs of Theorem~5}
\put(2.0,1.0){\small and Theorem~6: $n>3.$}
\end{picture}
\end{figure}

Let $\r$ be an optimal order for $\R$.  Contrasting two alternatives in the
manner described in the formulation of SCM will be called {\it
confrontation.} Then

(A) Confronting $X_1$ and $X_4$ and using SCM, we get $X_1\gro X_4.$

(B) Confronting $X_3$ and $X_4$ and assuming $X_4\gero X_3$, we get
$X_3\gro X_1$, in contradiction to (A). Therefore $X_3\gro X_4$.

(C) Confronting $X_2$ and $X_3$, we get $X_2\gro X_3.$

(D) Confronting $X_1$ and $X_2$ and using (B), we get $X_1\gro X_2.$

Thus, it follows from SCM that the restriction of any optimal order to
$\{X_1, X_2,$ $X_3, X_4\}$ is the transitive closure of $\{(X_1,X_2),
(X_2,X_3), (X_3,X_4)\}$.

Consider the weak order $\r'$ that is obtained from $\r$ by
interchanging $X_1$ and $X_2.$ Then for all $r_{ij}^p$ defined in
$\R$, $\sign(\r i-\r j) =\sign(\r'i-\r'j)$, and $\r'$ is optimal too,
since the operator is indifferent to the degree of resulting
preferences by our assumption. On the other hand, $r'$ violates SCM
(see (D)). This contradiction proves the desired statement.

2. $n=3.$ Consider the following $\R:$
$r_{12}^1=r_{13}^1=r_{\max},$
$r_{21}^1=r_{31}^1=r_{\min};$ all other $r_{ij}^p$ with $i\ne
j$ are undefined (Figure A4).

\begin{figure}[htb]
\setlength{\unitlength}{5mm}
\begin{picture}(26,12.6)
\put( 4, 5){\line(1,0){7}}
\put( 4,11){\line(1,0){7}}
\put( 4, 5){\line(0,1){6}}
\put(11, 5){\line(0,1){6}}
\put(4,5){\begin{picture}(7,6)
             \put(1.5,2.0){\circle{0.2}} 
             \put(3.5,4.0){\circle{0.2}} 
             \put(5.5,2.0){\circle{0.2}} 
             \put(1.1,1.0){$X_2$}
             \put(3.2,4.3){$X_1$}
             \put(5.1,1.0){$X_3$}
             \put(3.3,3.8) {\vector(-1,-1){1.6}} 
             \put(3.7,3.8) {\vector( 1,-1){1.6}} 
             \put(3.0,-1){$R^{(1)}$}
          \end{picture}}
\put(14, 5){\line(1,0){7}}
\put(14,11){\line(1,0){7}}
\put(14, 5){\line(0,1){6}}
\put(21, 5){\line(0,1){6}}
\put(14,5){\begin{picture}(7,6)
             \put(1.5,2.0){\circle{0.2}} 
             \put(3.5,4.0){\circle{0.2}} 
             \put(5.5,2.0){\circle{0.2}} 
             \put(1.1,1.0){$X_2$}
             \put(3.2,4.3){$X_1$}
             \put(5.1,1.0){$X_3$}
             \put(0.5,-1){$R^{(2)}=\cdots=R^{(m)}$}
          \end{picture}}
\put(2.0,2.0){\small{\sc Figure A4.}  $\R=(R^{(1)}\cdc R^{(m)})$
in the proofs of Theorem~5}
\put(2.0,1.0){\small and Theorem~6: $n=3.$}
\end{picture}
\end{figure}

Let $\r$ be an optimal order for $\R$.  Then using SCM and confronting
$X_1$ and $X_2$, we have $X_1\gro X_2;$ confronting $X_2$ and $X_3$, we
get $X_2\sro X_3.$ On the other hand, indifference to the degree of
resulting preferences implies that the orders determined by $X_1\gro
X_2\gro X_3$ and $X_1\gro X_3\gro X_2$ are optimal too, in
contradiction to SCM. The theorem is proved.  Finally, note that the
latter argument can be extended to the case $n>3.$ However, we have
preferred to give another proof for that case, since it demonstrates
that indifference to the degree of resulting preferences allows (in
some cases) to set $X_j\gro X_i$ whereas by SCM, $X_i\gro X_j$.
\bigskip

{\sc Proof of Corollary 1:} It suffices to show that these operators
are indifferent to the degree of resulting preferences. This is
obvious.
\bigskip

{\sc Proof of Theorem 6:} Assume that one of these operators satisfies
SCM.

1. Let $n>3$. Consider the following $\R$, the same as in the proof of
Theorem~5:
$r_{13}^1=r_{34}^1=r_{24}^1=r_{\max},$
$r_{31}^1=r_{43}^1=r_{42}^1=r_{\min};$ all other $r_{ij}^p$ with $i\ne
j$ are undefined (Figure~A3). Let $\r$ be an optimal order for
$\R$.  As have been shown in the proof of Theorem~5, SCM implies
that the restriction of any optimal order to
$\{X_1,X_2,X_3,X_4\}$ is the transitive closure of $\{(X_1,X_2),
(X_2,X_3), (X_3,X_4)\}$.

Consider the weak order $\r'$ that is obtained from $\r$ by
interchanging $X_1$ and $X_2.$ Let $f_k$ be the objective functions of
$\WQA(R,C'_k)$, $k\in\{4,5,6\}$ and let $f_{\rm NDW}$ be the objective
function of ``Net-Difference-WALB". Then for $\R$,
$$
f_5(\r)=(\r 1-\r 3+1)r_{\max}+(\r 2-\r 4+1)r_{\max}+(\r 3-\r 4
+1)r_{\max}
$$
$$
=(\r 1+\r 2-2\r 4+3)r_{\max}=f_5(\r');
$$
$$
f_6(\r)=(\r 3-\r 1-1)r_{\min}+(\r 4-\r 2-1)r_{\min}+(\r 4-\r 3
-1)r_{\min}
$$
$$
=(2\r 4-\r 1-\r 2-3)r_{\min}=f_6(\r');
$$
$f_4(\r)=f_4(\r')$ by Theorem~1, and
$$
f^{}_{\rm NDW}(\r)=0=f^{}_{\rm NDW}(\r').
$$

We see that all these operators equalize $\r'$ and $\r$ for $\R$, and
thus $\r'$ is also optimal, which contradicts SCM.

2. $n=3.$ Consider the following $\R$, the same as in the proof of
Theorem~5:
$r_{12}^1=r_{13}^1=r_{\max},$
$r_{21}^1=r_{31}^1=r_{\min};$ all other $r_{ij}^p$ with $i\ne
j$ are undefined (Figure~A4). By SCM, a unique optimal weak order
$\r$ is determined by $X_1\gro X_2\sro X_3.$ On the other hand,
each of the four operators under consideration equalizes $\r$ and
the orders determined by $X_1\gro X_2\gro X_3$ and $X_1\gro
X_3\gro X_2$. Hence they are also optimal, which contradicts SCM.
The theorem is proved.  Here the final remark in the proof of
Theorem~5 is applicable as well.
\bigskip

{\sc Proof of Theorem 7:} Let $n>4$. Consider the following $\R:$ for
all $(i,j)\in\{(k,\ell):\; k,\ell\in\{1\cdc n\},\; k<\ell\}\setminus
\{(1,2),(1,4),(2,3)\},$ $r_{ij}^1=r_{\max}$ and $r_{ji}^1=r_{\min};$
all other $r_{ij}^p$ with $i\ne j$ are undefined (Figure A5).

\begin{figure}[htb]
\setlength{\unitlength}{5mm}
\begin{picture}(26,25.6)
\put( 3, 6){\line(1,0){10}}
\put( 3,24){\line(1,0){10}}
\put( 3, 6){\line(0,1){18}}
\put(13, 6){\line(0,1){18}}
\put(3,6){\begin{picture}(10,18)
            \put(1.5,9.5){\begin{picture}(8,8)
               \put(1.5,6.0){\circle{0.2}} 
               \put(1.5,4.0){\circle{0.2}} 
               \put(5.5,5.0){\circle{0.2}} 
               \put(5.5,2.0){\circle{0.2}} 
               \put(1.2,6.3){$X_1$}
               \put(1.1,3.0){$X_3$}
               \put(5.2,5.3){$X_2$}
               \put(5.2,1.0){$X_4$}
               \put(1.5,5.7) {\vector(0,-1){1.4}} 
               \put(1.8,3.87){\vector(2,-1){3.3}} 
               \put(5.5,4.7) {\vector(0,-1){2.3}} 
               \put(3.5,3.75) {\oval(7,7.5)}
            \end{picture}}
            \put(5.0,7.0){\circle{0.2}} 
            \put(5.0,2.7){\circle{0.2}} 
            \put(4.7,7.3){$X_5$}
            \put(4.6,1.7){$X_n$}
            \put(5.0,6.7){\vector(0,-1){0.9}}
            \put(5.0,4.0){\vector(0,-1){1.0}}
            \put(4.83,6.66){\line(-1,-2){0.83}}
            \put(4.0,5.0){\vector( 1,-1){0.83}}
            \put(4.7,6.6){\line(-3,-4){1.2}}
            \put(3.5,5.0){\vector(2,-3){1.3}}
            \put(5.0,5.0){\circle*{0.15}}
            \put(5.0,5.5){\circle*{0.15}}
            \put(5.0,4.5){\circle*{0.15}}
            \put(5.0,4.75){\oval(4,7.5)}
            \linethickness{.5mm}
            \put(5.0,9.5) {\vector(0,-1){0.9}}
            \put(4.6,-1){$R^{(1)}$}
          \end{picture}}
\put(16, 6){\line(1,0){10}}
\put(16,24){\line(1,0){10}}
\put(16, 6){\line(0,1){18}}
\put(26, 6){\line(0,1){18}}
\put(16,6){\begin{picture}(10,18)
            \put(1.5,9.5){\begin{picture}(8,8)
               \put(1.5,6.0){\circle{0.2}} 
               \put(1.5,4.0){\circle{0.2}} 
               \put(5.5,5.0){\circle{0.2}} 
               \put(5.5,2.0){\circle{0.2}} 
               \put(1.2,6.3){$X_1$}
               \put(1.1,3.0){$X_3$}
               \put(5.2,5.3){$X_2$}
               \put(5.2,1.0){$X_4$}
            \end{picture}}
            \put(5.0,7.0){\circle{0.2}} 
            \put(5.0,2.7){\circle{0.2}} 
            \put(4.7,7.3){$X_5$}
            \put(4.6,1.7){$X_n$}
            \put(5.0,5.0){\circle*{0.15}}
            \put(5.0,5.5){\circle*{0.15}}
            \put(5.0,4.5){\circle*{0.15}}
            \put(2.1,-1){$R^{(2)}=\cdots=R^{(m)}$}
          \end{picture}}
\put(2.0,3.0){\small{\sc Figure A5.}  $\R=(R^{(1)}\cdc R^{(m)})$
in the proof of Theorem~7. The bold}
\put(2.0,2.0){\small arrow between two subsets designates that for
every $X_i$ in the first subset}
\put(2.0,1.0){\small and every $X_j$ in the second subset,
$r_{ij}^p=r_{\max}$ and $r_{ji}^p=r_{\min}.$}
\end{picture}
\end{figure}

Then $t_1=t_2=n-3,$ $t_3=n-4,$ $t_4=n-6,$ and for $i=5\cdc n,$
$t_i=n+1-2i.$

Let us prove that for $\R$ there exists only one weak order $\r$
satisfying SCM and that it is determined by $X_1\gro X_2\gro X_3\gro
X_4\gro X_5\gro\cdots\gro X_n.$ Indeed, we have the following.

(A) For any $i\in\{1\cdc 4\}$ and $j\in\{5\cdc n\}$, confronting $X_i$
and $X_j$ yields $X_i\gro X_j$.

(B) For any $i\in\{5\cdc n-1\}$ and $j\in\{i+1\cdc n\}$, confronting
$X_i$ and $X_j$ yields $X_i\gro X_j$.

(C) Confronting $X_1$ and $X_4$, we get $X_1\gro X_4$.

(D) Confronting $X_2$ and $X_3$, we get $X_2\gro X_3$.

(E) Assuming $X_4\gero X_3$ and confronting $X_3$ and $X_4$, we get
$X_3\gro X_1$ in contradiction to (C). Therefore $X_3\gro X_4$.

(F) Confronting $X_1$ and $X_2$ and using (E), we have $X_1\gro X_2$.

Now consider $\r'$ determined by $X_1\srop X_2\grop X_3\grop
X_4\grop\cdots\grop X_n$ and let $f(\cdot)$ be the objective function
of the $\beta$-LS method. It can be shown that for this $\R$
$$
f(\r)-f(\r')=4\beta^2(n-5).
$$
Therefore, $\r$ is not uniquely optimal for $r=5$ and is not optimal at
all for $n>5$. Hence the $\beta$-LS method violates Self-Consistency,
and the theorem is proved.
\bigskip

{\sc Proof of Theorem 8:} Suppose that the conditions of the nonstrict
part of Self-Consistent Monotonicity are satisfied but $X_j\gro X_i,$
where $\r$ is the optimal weak order for the generalized row sum
method.  Consider the $i$th and $j$th equations of (4):
\begin{eqnarray}{\jot.5ex}
x\_i&=
\suml_{r_{ik}^p\in\R_i'}\hspace{-4pt}\big((1+\ve mn)r_{ik}^p+\ve\cdot(x\_k-x\_i)\big)
+\suml_{r_{ik}^p\in\R_i''}\hspace{-4pt}\big((1+\ve mn)r_{\max}+
\ve\cdot(x\_k-x\_i)\big),
\nonumber\\
x\_j&=
\suml_{r_{jk}^p\in\R_j'}\hspace{-4pt}\big((1+\ve mn)r_{jk}^p+\ve\cdot(x\_k-x\_j)\big)
+\suml_{r_{jk}^p\in\R_j''}\hspace{-4pt}\big((1+\ve mn)r_{\min}+
\ve\cdot(x\_k-x\_j)\big).
\nonumber
\end{eqnarray}

For every $r_{ik}^p\in\R_i'$, by $r_{j\tilde k}^{\tilde
p}$ denote the corresponding comparison outcome in $\R_j'.$ After
subtraction, we get
\begin{eqnarray}
x\_i-x\_j&=&
\suml_{r_{ik}^p\in\R_i'}
\big(
(1+\ve mn)(r_{ik}^p-r_{j\tilde k}^{\tilde p})
+\ve\cdot(x\_k-x\_{\tilde k})
+\ve\cdot(x\_j-x\_i)
\big)
\nonumber\\
&+&\suml_{r_{ik}^p\in\R_i''}\big((1+\ve mn)r_{\max}+
\ve\cdot(x\_k-x\_i)\big)
\nonumber\\
&+&\suml_{r_{jk}^p\in\R_j''}\big((1+\ve mn)r_{\max}+
\ve\cdot(x\_j-x\_k)\big).
\nonumber
\end{eqnarray}

Suppose that $\ve$ is reasonable and positive. By our
assumptions, all terms in the right-hand side are non-negative,
whereas the left-hand side is negative. This contradiction proves
that $X_i\gero X_j.$ The strict part of Self-Consistent
Monotonicity and Self-Consistency can be proved in the same way.
\bigskip

{\it Acknowledgements.} This work was supported by the Russian
Foundation for Basic Research. Partial research support from the
European Community under Grant No.~ACE-91-R02 is also gratefully
acknowledged.

\bigskip\bigskip
\centerline{REFERENCES}
\bigskip

\lfi{\sc Arditti, D.} (1989): ``Un Nouvel Algorithme de Recherche d'un
Ordre Induit par des Comparaisons par Paires,'' in {\it Troisi\`emes
Journ\'ees Internationales Analyse des Donn\'ees et Informatique},
Versailles, 4--7 oct., 1983, T.1, ed. E.~Diday et al. Le Chensay
e.a., s.a., 272--297.

\lfi{\sc Barth\'elemy, J.~P.} (1989): ``Social Welfare and Aggregation
Procedures: Com\-binatorial and Algorithmic Aspects,'' in {\it
Application of Combinatorics and Graph Theory to the Biological and
Social Sciences}, ed. by F.~Roberts. New York: Springer-Verlag, 39--73.

\lfi{\sc Belkin, A.~R., and M.~S.~Levin} (1990): {\it Decision Making:
Combinatorial Models of Information Approximation} [in Russian].
Moscow: Nauka.

\lfi{\sc Chebotarev, P.~Yu.} (1988): ``Two Methods of Ranking Objects
on the Basis of an Arbitrary Set of Paired Comparisons'' [in Russian].
VINITI manuscript No. 5879-B88, Moscow.

\lfi\sameauthor (1989): ``Generalization of the Row Sum Method for
Incomplete Paired Comparisons,'' {\it Automation and Remote Control},
50, 1103--1113.

\lfi\sameauthor (1990): ``On Some Optimization Methods for Aggregating
Preferences,'' in {\it Problems of Computerization and Statistics in
Applied Sciences} [in Russian]. Moscow: VNIISI Akad. Nauk SSSR, 67--72.

\lfi\sameauthor (1994): ``Aggregation of Preferences by the Generalized
Row Sum Method,'' {\it Mathematical Social Sciences}, 27, 293--320.

\lfi{\sc Chebotarev, P.~Yu., and E. V. Shamis} (1994):
``Characteristic Conditions for Aggregating Incomplete Preferences,''
in {\it Proceedings of the Fifth International Conference on
Statistical and Discrete Data Analysis, Odessa, September 1994} [in
Russian] (in press).

\lfi\sameauthor (1996):  ``Preference Fusion When the Number of
Alternatives Exceeds Two: Indirect Scoring Procedures,'' in {\it
Proceedings of Workshop on Foundations of Information/Decision Fusion:
Application to Engineering Problems, August 7-9, 1996, Washington
D.C.}. Lafayette, LA: Acadiana Printing, Inc., 20--32.

\lfi{\sc Cook, W.~D., and M.~Kress} (1992): {\it Ordinal Information
and Preference Structures: Decision Models and Applications}. Englewood
Cliffs, New Jersey: Prentice-Hall.

\lfi{\sc Critchlow, D.~E.} (1985): {\it Metric Methods for Analyzing
Partially Ranked Data}. Berlin-Heidelberg: Springer-Verlag.

\lfi{\sc Crow, E.~L.} (1990): ``Ranking Paired Contestants,'' {\it
Communications in Statistics. Simulation and Computation}, 19,
749--769.

\lfi\sameauthor (1993): ``Ranking from Paired Comparisons by Minimizing
Inconsistency,'' in {\it Probability Models and Statistical Analyses
for Ranking Data}, ed. by M.~A.~Fligner, and J.~S.~Verducci.  New York:
Springer-Verlag, 289--293.

\lfi{\sc David, H.~A.} (1987): ``Ranking from Unbalanced
Paired-Comparison Data,'' {\it Bio\-metrika}, 74, 432--436.

\lfi\sameauthor (1988): {\it The Method of Paired Comparisons.}
2nd ed. London: Griffin.

\lfi{\sc Frey, J.-J., and A.~Yehia-Alcoutlabi} (1986): ``Comparaisons
par Paires: Une Interpr\'etation et Une G\'en\'eralisation de la
M\'ethode des Scores,'' {\it RAIRO Recherche Operationnelle}, 20,
213--227.

\lfi{\sc Hubert, L.} (1976): ``Seriation Using Asymmetric Proximity
Measures,'' {\it British Journal of Mathematical and Statistical
Psychology}, 29, 32--52.

\lfi{\sc Kano, M., and A.~Sakamoto} (1983): ``Ranking the Vertices of a
Weighted Digraph Using the Length of Forward Arcs,'' {\it
Networks}, 13, 143--151.

\lfi\sameauthor (1985): ``Ranking the Vertices of a Paired Comparison
Digraph,'' {\it SIAM Journal on Algebraic and Discrete Methods}, 6,
79--92.

\lfi{\sc Kemeny, J.} (1959): ``Mathematics without Numbers,'' {\it
Daedalus}, 88, 571--591.

\lfi{\sc Kendall, M. G.} (1955): ``Further Contributions to the Theory
of Paired Comparisons,'' {\it Bio\-metrika}, 11, 43--62.

\lfi\sameauthor (1970): {\it Rank Correlation Methods}. 4th ed.
London: Griffin.

\lfi{\sc Litvak, B.~G.} (1982): {\it Expert Information: Methods of
Gathering and Analysis} [in Russian]. Moscow: Radio i Svyaz.

\lfi{\sc Slater, P.} (1961): ``Inconsistencies in a Schedule of Paired
Comparisons,'' {\it Bio\-metrika}, 48, 303--312.

\lfi{\sc Thompson, M.} (1975): ``On any Given Sunday: Fair Competitor
Orderings with Maximum Likelihood Methods,'' {\it Journal of the
American Statistical Association}, 70, 536--541.

\lfi{\sc Van Blokland-Vogelesang, R.} (1991): {\it Unfolding and
Group Consensus Ranking for Individual Preferences.} Leiden: DSWO
Press.

\lfi{\sc Weiss, H.~J., and J.~Y.~Assous} (1987): ``Reduction in Problem
Size for Ranking Alternatives in Group Decision-Making,'' {\it
Computers \& Operations Research}, 14, 55--65.

\lfi{\sc Young, H.~P.} (1986): ``Optimal Ranking and Choice from
Pairwise Comparisons,'' in {\it Information Pooling and Group Decision
Making}, ed. by B.~Grofman, and G.~Owen. Greenwich, Connecticut: JAI
Press, 113--122.

\end{document}